\documentclass[12pt,draft]{amsart}
\usepackage{amsmath,amsthm,latexsym,amscd,amsbsy,amssymb}
\chardef\bslash=`\\ 

\makeatletter
\def\verbatim{\interlinepenalty\@M \@verbatim
  \leftskip\@totalleftmargin\advance\leftskip2pc
  \frenchspacing\@vobeyspaces \@xverbatim}
\makeatother
\newcounter{rmnum}

\newtheorem{thm}{Theorem}[section]
\newtheorem{cor}[thm]{Corollary}

\newtheorem{pro}[thm]{Proposition}

\theoremstyle{definition}

\theoremstyle{remark}

\numberwithin{equation}{section}

\begin{document}


\title[Complemented subspaces of locally convex direct sums of Banach spaces]
{Complemented subspaces of locally convex direct sums of Banach spaces}
\author{Alex Chigogidze}
\address{Department of Mathematics and Statistics,
University of Saskatche\-wan,
McLean Hall, 106 Wiggins Road, Saskatoon, SK, S7N 5E6,
Canada}
\email{chigogid@math.usask.ca}
\thanks{Author was partially supported by NSERC research grant.}



\begin{abstract}{We show that a complemented subspace of a locally
convex direct sum of an uncountable collection of Banach spaces
is a locally convex direct sum of complemented subspaces
of countable subsums. As a corollary we prove that a complemented
subspace of a locally convex direct sum of arbitrary collection of
$\ell_{1}(\Gamma )$-spaces is isomorphic to a locally convex direct
sum of $\ell_{1}(\Gamma )$-spaces.} 
\end{abstract}

\maketitle
\markboth{A.~Chigogidze}{Complemented subspaces of locally
convex direct sums of Banach spaces}

\section{Introduction}\label{S:intro}
In 1960 A.~Pelczynski proved \cite{pel60} that complemented subspaces
of $\ell_{1}$ are isomorphic to $\ell_{1}$. In \cite{kot66} G.~K\"{o}the
generalized this result to the non-separable case. Later, while answering K\"{o}the's
question about precise description of projective spaces in the category of
(LB)-spaces, P.~Doma\'{n}ski showed \cite{dom92} that complemented
subspaces of locally convex direct sums of countable collections of
$\ell_{1}(\Gamma )$-spaces have the same structure, i.e. are
isomorphic to locally convex direct sums of countable
collections of $\ell_{1}(\Gamma )$-spaces.

Below we complete this series of statements by showing
(Corollary \ref{C:ellspaces}) that countability assumption in Doma\'{n}ski's
result is not essential. More precisely, we prove 
that complemented subspaces of a locally convex direct sums of
arbitrary collections of $\ell_{1}(\Gamma )$-spaces are isomorphic to
locally convex direct sums of $\ell_{1}(\Gamma )$-spaces. This is
obtained as a corollary of our main result (Theorem \ref{T:complementedsub})
stating that complemented subspaces of locally convex direct
sums of arbitrary collections of Banach spaces are isomorphic to
locally convex direct sums of complemented subspaces of countable subsums.


\section{Results}
Below we work with locally convex
direct sums $\bigoplus\{ B_{t} \colon t \in T\}$ of
uncountable collections of Banach spaces $B_{t}$, $t \in T$.
Recall that if $S \subseteq R \subseteq T$,
then $\bigoplus\{ B_{t} \colon t \in S\}$ can be
canonically identified with the subspace
\[\left\{ \{ x_{t} \colon t \in R\} \in
\bigoplus\{ B_{t} \colon t \in R\} \colon x_{t} = 0\;
\text{for each}\; t \in R-S\right\} \]

\noindent of
$\bigoplus\{ B_{t} \colon t \in R\}$. The corresponding
inclusion is denoted by $i_{S}^{R}$. The following
statement is used in the proof of
Theorem \ref{T:complementedsub}.

\begin{pro}\label{P:injspectral}
Let $r \colon \bigoplus\{ B_{t} \colon t \in T\}
\to \bigoplus\{ B_{t} \colon t \in T\}$ be a
continuous linear map of a locally convex direct sum of an
uncountable collection of
Banach spaces into itself. Let also $A$ be a countable subset of $T$.
Then there exists a countable subset $S \subseteq T$ such
that $A \subseteq S$
and
$r\left(\bigoplus\{ B_{t} \colon t \in S\}\right)
\subseteq \bigoplus\{ B_{t} \colon t \in S\}$.
\end{pro}
\begin{proof}
Let $\exp_{\omega}T$ denote the set of all countable subsets
of the indexing set $T$.
Consider the following relation
\[
{\mathcal L} = \left\{ (A,C) \in (exp_{\omega}T)^2
\colon A \subseteq C \; \text{and}\; r\left(\bigoplus\{
B_{t} \colon t \in A \}\right) \subseteq
\bigoplus\{ B_{t} \colon t \in C\}\right\} .\]

\noindent We need to verify the following three properties of the above
defined relation.
\medskip

{\em Existence}. If $A \in \exp_{\omega}T$, then there
exists $C \in \exp_{\omega}T$
such that $(A,C) \in {\mathcal L}$.
\medskip

{\em Proof}. First of all let us make the following observation.
\medskip

{\bf Claim.} {\em For each $j \in T$ there exists a
finite subset $C_{j} \subseteq T$ such that
$r\left( B_{j}\right) \subseteq
\bigoplus\{ B_{t} \colon t \in C_{j}\}$.}
\medskip

{\em Proof of Claim}. The unit ball $K = \{ x \in B_{j} \colon
||x||_{j} \leq 1\}$ (here $||\cdot ||_{j}$ denotes the
norm of the Banach space $B_{j}$) being bounded in $B_{j}$
is, by \cite[Theorem 6.3]{scha66}, bounded in
$\bigoplus\{ B_{t} \colon t \in T\}$. Continuity of $r$
guarantees that $r(K)$ is also bounded in
$\bigoplus\{ B_{t} \colon t \in T\}$.
Applying \cite[Theorem 6.3]{scha66} once again, we
conclude that there exists
a finite subset $C_{j} \subseteq T$ such that
$r(K) \subseteq \bigoplus\{ B_{t} \colon t \in C_{j}\}$.
Finally the linearity of $r$ implies that
$r\left( B_{j}\right) \subseteq
\bigoplus\{ B_{t} \colon t \in C_{j}\}$ and proves the Claim.

Let now $A \in \exp_{\omega}T$. For each $j \in A$,
according to Claim,
there exists a finite subset $C_{j} \subseteq T$ such that
$r\left( B_{j}\right) \subseteq
\bigoplus\{ B_{t} \colon t \in C_{j}\}$. Without loss of
generality we may assume that $A \subseteq C_{j}$ for each $j \in A$.
Let $C = \cup\{ C_{j} \colon j \in A\}$.
Clearly $C$ is countable, $A \subseteq C$ and
$r\left( B_{j}\right) \subseteq \bigoplus\{ B_{t}
\colon t \in C\}$ for each $j \in A$. This
guarantees that $r\left( \bigoplus\{ B_{t} \colon t \in A \}\right)
\subseteq \bigoplus\{ B_{t} \colon t \in C\}$ and shows
that $(A,C) \in {\mathcal L}$.
\medskip

{\em Majorantness}. If $(A,C) \in {\mathcal L}$,
$D \in \exp_{\omega}T$ and $C \subseteq D$, then $(A,D) \in {\mathcal L}$.
\medskip

{\em Proof}. Condition $(A,C) \in {\mathcal L}$ implies that
$r\left(\bigoplus\{ B_{t} \colon t \in A\}\right)
\subseteq \bigoplus\{ B_{t} \colon t \in C\}$. The inclusion
$C \subseteq D$ implies that
$\bigoplus\{ B_{t} \colon t \in C\}
\subseteq \bigoplus\{ B_{t} \colon t \in D\}$. Consequently
$r\left(\bigoplus\{ B_{t} \colon t \in A\}\right)
\subseteq \bigoplus\{ B_{t} \colon t \in C\} \subseteq
\bigoplus\{ B_{t} \colon t \in D\}$,
which means that $(A,D) \in {\mathcal L}$.
\medskip

{\em $\omega$-closeness}. Suppose that $(A_{i},C) \in
{\mathcal L}$ and $A_{i} \subseteq A_{i+1}$ for each
$i \in \omega$. Then $(A,C) \in {\mathcal L}$, where
$A = \cup \{ A_i \colon i \in \omega \}$. 
\medskip

{\em Proof}. Consider the following inductive sequence

\[ \bigoplus\{ B_{t} \colon t \in A_{0}\}
\xrightarrow{i_{A_{0}}^{A_{1}}}\cdots\xrightarrow{ }
\bigoplus\{ B_{t} \colon t \in A_{i}\}
\xrightarrow{i_{A_{i}}^{A_{i+1}}}\bigoplus\{ B_{t}
\colon t \in A_{i+1}\} \xrightarrow{ }\cdots\]

\noindent limit of which is isomorphic to
$\bigoplus\{ B_{t} \colon t \in A\}$ (horizontal arrows
represent canonical inclusions). Since $r\left(\bigoplus\{ B_{t} \colon t
\in A_{i}\}\right) \subseteq \bigoplus\{ B_{t}
\colon t \in C\}$ for each $i \in \omega$ (assumption
$(A_{i},C) \in {\mathcal L}$), it follows that

\begin{multline*}
r\left( \bigoplus\{ B_{t} \colon t \in A\}\right) =
r\left(  \injlim\left\{ \bigoplus\{ B_{t} \colon t
\in A_{i}\} ,  i_{A_{i}}^{A_{i+1}},
\in \omega\right\}\right) \subseteq\\
 \bigoplus\{ B_{t} \colon t \in C\} .
\end{multline*}

\noindent This obviously means that $(A,C) \in {\mathcal L}$ as required.
\medskip

According to \cite[Proposition 1.1.29]{chibook96} the set of
${\mathcal L}$-reflexive elements of $\exp_{\omega}T$ is
cofinal in $\exp_{\omega}T$. An element $S \in \exp_{\omega}T$
is ${\mathcal L}$-reflexive if $(S,S) \in {\mathcal L}$.
In our situation
this means that the given countable subset $A$ of $T$ is
contained in a larger countable subset $S$ such that
$r\left(\bigoplus\{ B_{t} \colon t \in S\}\right)
\subseteq \bigoplus\{ B_{t} \colon t \in S\}$. Proof is completed.
\end{proof}


\begin{thm}\label{T:complementedsub}
Let $T$ be an uncountable set. A complemented subspace of a locally
convex direct sum $\bigoplus\{ B_{t} \colon t \in T\}$ of
Banach spaces $B_{t}$, $t \in T$, is
isomorphic to a locally convex direct sum
$\bigoplus\{ F_{j} \colon j \in J\}$,
where $F_{j}$ is a complemented subspace of the countable sum
$\bigoplus\{ B_{t} \colon t \in T_{j}\}$
where $|T_{j}| = \omega$
for each $j \in J$.
\end{thm}
\begin{proof}
Let $X$ be a complemented subspace of the sum
$B = \bigoplus\{ B_{t} \colon t \in T\}$.
Choose a continuous
linear map $r \colon B \to X$ such that $r(x) = x$ for each
$x \in X$. Let us agree that a subset $S \subseteq T$ is called
$r$-admissible if
$r\left(\bigoplus\{ B_{t} \colon t \in S\}\right)
\subseteq \bigoplus\{ B_{t} \colon t \in S\}$.

For a subset $S \subseteq T$,
let $\displaystyle X_{S} = r\left(\bigoplus
\{ B_{t} \colon t \in S\}\right)$.

{\bf Claim 1}.
{\em If $S \subseteq T$ is an $r$-admissible, then
$\displaystyle X_{S} = X \bigcap \left(\bigoplus\{ B_{t}
\colon t \in S\}\right)$.}
\medskip

{\em Proof}. Indeed, if $y \in X_{S}$, then there exists a point
$\displaystyle x \in \bigoplus\{ B_{t} \colon t \in S\}$ such that
$r(x) = y$. Since $S$ is $r$-admissible, it follows that
\[\displaystyle y = r(x) \in r\left(\bigoplus
\{ B_{t} \colon t \in S\}\right) \subseteq
\bigoplus\{ B_{t} \colon t \in S\} .\]
\noindent Clearly, $y \in X$. This shows that
$\displaystyle X_{S} \subseteq X \bigcap \left(\bigoplus\{ B_{t}
\colon t \in S\}\right)$.

Conversely, if
$\displaystyle y \in X \bigcap \left(\bigoplus\{ B_{t}
\colon t \in S\}\right)$, then $y \in X$ and hence,
by the property of $r$,
$y = r(y)$. Since $\displaystyle y \in \bigoplus\{ B_{t}
\colon t \in S\}$, it follows that
$\displaystyle y = r(y) \in r\left(\bigoplus\{ B_{t}
\colon t \in S\}\right) = X_{S}$.
\medskip

{\bf Claim 2}.
{\em The union of an arbitrary collection of $r$-admissible
subsets of $T$ is $r$-admissible.}
\medskip

{\em Proof}. Straightforward
verification based of the definition of the $r$-ad\-missi\-bi\-lity.
\medskip

{\bf Claim 3.}
{\em Every countable subset of $T$ is contained in a countable
$r$-admissible subset
of $T$.}

{\em Proof}. This follows from Proposition \ref{P:injspectral}
applied to the map $r$.
\medskip

{\bf Claim 4.}
{\em If $S\subseteq T$ is an $r$-admissible subset of $T$, then
$r_{S}(x) = x$ for each point $x \in X_{S}$, where 
$\displaystyle r_{S} = r\left|\left( \bigoplus\{
B_{t} \colon t \in S\}\right.\right) \colon \bigoplus\{
B_{t} \colon t \in S\} \to X_{S}$.}

{\em Proof}. This follows from the corresponding property of the map $r$.
\medskip

Before we state the next property of $r$-admissible sets
note that if
$S \subseteq R \subseteq T$, then the map
\[ \pi_{S}^{R} \colon \bigoplus\{ B_{t} \colon t \in R\} \to
\bigoplus\{ B_{t} \colon t \in S\} ,\]
defined by letting
\[
\pi_{S}^{R}\left(\{ x_{t} \colon t \in R\}\right) = 
\begin{cases}
x_{t} ,\;\text{if}\; t \in S\\
0\; , \;\;\text{if}\; t \in R-S ,
\end{cases}
\]
\noindent is continuous and linear.
\medskip

{\bf Claim 5.} {\em Let $S$ and $R$ are $r$-admissible subsets
of $T$ and $S \subseteq R$.
Then $X_{S}$ is a complemented subspace in $X_{R}$ and $X_{R}/_{X_{S}}$
is a complemented subspace in $\bigoplus\{ B_{t} \colon t \in R-S\}$.}
\medskip

{\em Proof}. Consider the following commutative diagram

\[
\begin{CD}
\bigoplus\{ B_{t} \colon t \in R\} @>r_{R}>> X_{R}\\
@V\pi^{R}_{R-S}VV @VVpV\\
\oplus\{ B_{t} \colon t \in R-S\} =\oplus\{ B_{t}
\colon t \in R\} /_{\oplus\{ B_{t} \colon t \in S\}} @>q>> X_{R}/_{X_{S}}
\end{CD}
\]


\noindent in which $p$ is the canonical map and $q$ is defined
on cosets by letting (recall that $r_{R}\left( \bigoplus\{
B_{t} \colon t \in S \}\right) = X_{S}$)
\[ q\left( x+\bigoplus\{ B_{t} \colon t \in S\}\right) =
r_{R}(x) +X_{S} \;\;\text{for each}\;\; x \in \bigoplus\{
B_{t} \colon t \in R\} .\]

Let us denote by $i_{R} \colon X_{R} \hookrightarrow \bigoplus\{
B_{t} \colon t \in R\}$
the natural inclusion and consider a map 
\[ j \colon X_{R}/_{X_{S}} \to \oplus\{ B_{t} \colon t \in R\}
/_{\oplus\{ B_{t} \colon t \in S\}}\]
defined by letting (in terms of cosets)

\[ j\left( x+X_{S}\right) = i(x) +\bigoplus\{ B_{t} \colon t \in S\}
\;\;\text{for each}\;\;  x\in X_{R} .\]

\noindent Note that $q\circ j = \operatorname{id}_{X_{R}/_{X_{S}}}$
(this follows from
the equality $r_{R}\circ i = \operatorname{id}_{X_{R}}$). In
particular, this shows that $X_{R}/_{X_{S}}$
is isomorphic to a complemented subspace of
$\bigoplus\{ B_{t} \colon t \in R-S\}$.

Finally consider the composition
$r_{R} \circ i^{R}_{R-S} \circ j \colon X_{R}/_{X_{S}} \to X_{R}$ and note that
\[ p\circ (r_{R} \circ i^{R}_{R-S} \circ j) =
p\circ r_{R} \circ i^{R}_{R-S} \circ j =
q \circ \pi_{R-S}^{R} \circ i^{R}_{R-S} \circ j =
q \circ \operatorname{id} \circ j = \operatorname{id}_{X_{R}/_{X_{S}}} .\]
This shows that $X_{S}$ is a complemented subspace of $X_{R}$ and
completes the proof of Claim 5.
\medskip

Let $|T| = \tau$. Then we can write
$T = \{ t_{\alpha} \colon \alpha < \tau\}$.
Since the collection of countable $r$-admissible
subsets of $T$ is cofinal in
$\exp_{\omega}T$ (see Claim 3), each element
$t_{\alpha} \in T$ is contained in a countable $r$-admissible
subset $A_{\alpha} \subseteq T$. According to Claim 2, the set
$T_{\alpha} = \bigcup\{ A_{\beta} \colon
\beta \leq \alpha\}$ is $r$-admissible for each $\alpha < \tau$.
Consider the inductive system ${\mathcal S} =
\{ X_{\alpha}, i_{\alpha}^{\alpha +1}, \tau\}$,
where $X_{\alpha} = X_{T_{\alpha}} = X\bigcap r\left(\bigoplus\{ B_{t}
\colon t \in T_{\alpha}\}\right)$ (see Claim 1)
and $i_{\alpha}^{\alpha +1} \colon X_{\alpha} \to X_{\alpha +1}$
denotes the natural inclusion for each $\alpha < \tau$.
For a limit ordinal number $\beta < \tau$ the space
$X_{\beta}$ is
isomorphic to the limit space of the direct system
$\left\{ X_{\alpha},
i_{\alpha}^{\alpha +1}, \alpha < \beta \right\}$
(verification of this fact is based on Claim 4 coupled with
the fact that $\bigoplus\{ B_{t} \colon t \in T_{\beta}\}$ is
isomorphic to the limit of the direct system
$\left\{ \bigoplus\{ B_{t} \colon t \in T_{\alpha}\} ,
i_{T_{\alpha}}^{T_{\alpha +1}}, \alpha < \beta\right\}$).
In particular, $X$ is isomorphic to the limit
of the inductive system
$\{ X_{\alpha}, i_{\alpha}^{\alpha +1}, \alpha <\tau\}$.

For each $\alpha < \tau$, according to Claim 5, the inclusion
$i_{\alpha}^{\alpha +1} \colon X_{\alpha} \to X_{\alpha +1}$ is
isomorphic to the inclusion
$X_{\alpha} \hookrightarrow X_{\alpha} \bigoplus
X_{\alpha +1}/_{X_{\alpha}}$. In this situation the
straightforward transfinite induction shows that $X$ is
isomorphic to the locally convex direct sum 
$X_{0}\bigoplus \left(\bigoplus \left\{ X_{\alpha +1}/_{X_{\alpha}} \colon
\alpha < \tau \right\}\right)$.

By construction, the set $T_{0}$ is countable and $X_{0}$ is a
complemented subspace of $\bigoplus\{ B_{t} \colon t \in T_{0}\}$.
Note also that for each $\alpha < \tau$ the
set $T_{\alpha +1}-T_{\alpha} = A_{\alpha +1}$ is countable
and $X_{\alpha +1}/_{X_{\alpha}}$ is a complemented subspace
of $\bigoplus\{ B_{t} \colon t \in A_{\alpha +1}\}$. This completes
the proof
of Theorem \ref{T:complementedsub}.
\end{proof}

The following statement, as was noted in the Introduction,
provides a complete description of complemented subspaces of
locally convex direct sums of uncountable collections of
$\ell_{1} (\Gamma )$-spaces.

\begin{cor}\label{C:ellspaces}
Let $X$ be a complemented subspace of
$\bigoplus\{ \ell_{1}(\Gamma_{t}) \colon t \in T\}$.
Then $X$ is isomorphic to $\bigoplus\{ \ell_{1}
(\Lambda_{i}) \colon i \in I\}$.
\end{cor}
\begin{proof}
For countable $T$ results follows from \cite{kot66}
and \cite{dom92}.
Let now $T$ is uncountable and $X$ be a complemented
subspace of
a locally convex direct sum
$\bigoplus\{ \ell_{1} (\Gamma_{t}) \colon t \in T\}$.
By Theorem \ref{T:complementedsub}, $X$ is isomorphic to
a locally convex direct sum
$\bigoplus\{ F_{j} \colon j \in J\}$,
where $F_{j}$ is a complemented subspace of the countable sum
$\bigoplus\{ \ell_{1}(\Gamma_{t}) \colon t \in T_{j}\}$
where $|T_{j}| = \omega$
for each $j \in J$. According to \cite{dom92},
$F_{j} = \bigoplus\{ \ell_{1}(\Lambda_{t}) \colon  \in t \in T_{j}\}$
for each $j \in J$. Consequently, $X$ is isomorphic to the
locally convex direct sum
$\bigoplus\{ \bigoplus\{ \ell_{1}(\Lambda_{t} \colon
t \in T_{j}\} \colon j \in J \}
= \bigoplus\{ \ell_{1}(\Lambda_{t}) \colon t \in
\cup\{ T_{j} \colon j \in J\} \}$ as required.
\end{proof}


\providecommand{\bysame}{\leavevmode\hbox to3em{\hrulefill}\thinspace}



\end{document}